\documentclass[11pt]{article}\textwidth 160mm\textheight 235mm
\usepackage{graphicx}
\usepackage{epstopdf}
\usepackage{amssymb,amsmath}
\usepackage{amsfonts}
\usepackage{tikz}
\usepackage{amsthm}
\usepackage{color}

\def\hat{\widehat}
\def\tilde{\widetilde}

\def\dom{{\rm dom}\,}

\def\Lra{\Longrightarrow}

\def\N{{\cal N}}

\def\O{{\cal O}}

\def\R{I\!\!R}
\def\N{I\!\!N}
\def\B{\mathbb B}

\def\ox{\overline{x}}
\def\oy{\overline{y}}

\def\disp{\displaystyle}

\def\Limsup{\mathop{{\rm Lim}\,{\rm sup}}}

\def\tto{\;{\lower 1pt \hbox{$\rightarrow$}}\kern -10pt
\hbox{\raise 2pt \hbox{$\rightarrow$}}\;}
\def\Hat{\widehat}
\def\Tilde{\widetilde}
\def\Bar{\overline}
\def\ra{\rangle}
\def\la{\langle}

\def\B{I\!\!B}
\def\h{\hfill\Box}
\def\R{\mathbb{R}}
\def\N{I\!\!N}
\def\ox{\bar{x}}
\def\oy{\bar{y}}

\def\gph{\mbox{\rm gph}\,}

\def\dom{\mbox{\rm dom}\,}

\def\lip{\mbox{\rm lip}\,}
\def\ssubreg{\mbox{\rm ssubreg}\,}
\def\subreg{\mbox{\rm subreg}\,}

\def\h{\hfill\triangle}
\def\dn{\downarrow}
\def\O{\Omega}
\def\ph{\varphi}

\def\st{\stackrel}
\def\oR{\Bar{\R}}
\def\lm{\lambda}
\def\gg{\gamma}

\def\al{\alpha}

\def\N{I\!\!N}

\def\sce{\setcounter{equation}{0}}
\begin{document}
\begin{center}
{\bf HIGHER-ORDER METRIC SUBREGULARITY AND ITS APPLICATIONS}\footnote{This research was partly supported by the National Science Foundation under grant DMS-12092508.}\\[2ex]
BORIS S. MORDUKHOVICH\footnote{Department of Mathematics, Wayne State University, Detroit, MI 48202 (boris@math.wayne.edu).} and WEI OUYANG\footnote{Department of Mathematics, Wayne State University, Detroit, MI 48202 (wei@wayne.edu).}
\end{center}
\small{\bf Abstract.} This paper is devoted to the study of metric subregularity and strong subregularity of any positive order $q$ for set-valued mappings in finite and infinite dimensions. While these notions have been studied and applied earlier for $q=1$ and---to a much lesser extent---for $q\in(0,1)$, no results are available for the case $q>1$. We derive characterizations of these notions for subgradient mappings, develop their sensitivity analysis under small perturbations, and provide applications to the convergence rate of Newton-type methods for solving generalized equations.\\[1ex]

\newtheorem{Theorem}{Theorem}[section]
\newtheorem{Proposition}[Theorem]{Proposition}
\newtheorem{Remark}[Theorem]{Remark}
\newtheorem{Lemma}[Theorem]{Lemma}
\newtheorem{Corollary}[Theorem]{Corollary}
\newtheorem{Definition}[Theorem]{Definition}
\newtheorem{Example}[Theorem]{Example}
\renewcommand{\theequation}{{\thesection}.\arabic{equation}}
\renewcommand{\thefootnote}{\fnsymbol{footnote}}

\normalsize
\section{Introduction}\sce

This paper mainly concerns the study and some applications of the notions of {\em higher-order} metric {\em sub}regularity and its {\em strong} subregularity counterpart. For definiteness, we use the number $q>0$ to indicate the {\em order/rate} of the corresponding regularity under consideration. Recall first that a set-valued mapping $F:X\rightrightarrows Y$ between Banach spaces is {\em metrically $q$-regular} at (better {\em around}) $(\bar x,\bar y)\in\gph F$ if there exist a number $\eta>0$ and neighborhoods $U$ of $\bar x$ and $V$ of $\bar y$ such that
\begin{equation}\label{1.1}
d\big(x;F^{-1}(y)\big)\le\eta\,d^q\big(y;F(x)\big)\;\mbox{ for all }\;x\in U\;\text{ and }\;y\in V,
\end{equation}
where $d(\cdot;\Omega)$ is the distance function associated with $\Omega$. It has been well recognized in nonlinear and variational analysis that metric regularity ($q=1$) and the equivalent notions of linear openness and Lipschitzian stability play an important role in optimization, control, equilibria, and various applications as documented, e.g., in the books \cite{bz05,dr09,bm06,rw} with many references therein. On the other hand, metric regularity often fails for broad classes of parametric variational systems given by the generalized equations in the sense of Robinson \cite{r}:
\begin{equation}\label{ge}
0\in f(x,y)+Q(y),
\end{equation}
where $f$ is single-valued while $Q(y)=\partial\varphi(y)$ is a set-valued mapping of the subdifferential/normal cone type generated by nonsmooth functions; see \cite{bm08} and also \cite{am,agh,bbem10,gmn09,u10} for more details and further results in this direction. However, this phenomenon does not appear if metric regularity is replaced by a weaker property of {\em metric subregularity} of $F$ at $(\bar x,\bar y)$ defined by
\begin{equation}\label{1.2}
d\big(x;F^{-1}(\bar y)\big)\le\eta\,d\big(\bar y;F(x)\big)\;\text{ for all }\;x\in U.
\end{equation}
Considering $\|x-\bar x\|$ instead of $d(x;F^{-1}(\bar y))$ in (\ref{1.2}), we get the notion of {\em strong metric subregularity}. In the aforementioned books and in an increasing number of papers, the reader can find more information about these subregularity properties, their {\em calmness} (resp.\  {\em isolated calmness}) equivalents for inverse mappings, as well as their various applications to optimization.

In \cite{bzh88,fq12,yyk, zn}, the authors studied the notion of {\em H\"older metric regularity}, which corresponds to (\ref{1.1}) with the replacement of $d(y;F(x))$ by $d^q(y;F(x))$ as $0<q<1$. Replacing $d(\bar y;F(x))$ by $d^q(\bar y;F(x))$ in (\ref{1.2}) as $0<q<1$ gives us the notion of {\em H\"older metric subregularity} considered recently in \cite{ggj,kkk,lm} from different viewpoints while without its {\em strong} counterpart.

It is essential to mention that there is no sense to study metric $q$-regularity of single-valued or set-valued mappings for $q>1$, since only constant mappings satisfy this property. However, it is not the case for $q$-subregularity that is equally important whenever $q>0$ as demonstrated in this paper, where---to the best of our knowledge---the notion of $q$-subregularity for $q>1$ is studied and applied for the first time in the literature.

In what follows we investigate both notions of metric $q$-subregularity and strong metric $q$-subregularity for any positive $q$ concentrating mainly on the higher-order case of $q>1$. In this way we derive verifiable sufficient conditions and necessary conditions for these notions of $q$-subregularity in terms of appropriate generalized differential constructions of variational analysis, study their behavior with respect to perturbations, and obtain their applications to the rate of convergence of Newton's and quasi-Newton methods for solving generalized equations.\vspace*{0.05in}

Accordingly, we organize the rest of the paper. Section~2 contains some preliminaries from variational analysis and generalized differentiation widely used in the formulations and proofs of the main results given below. Section~3 is devoted to a detailed study of $q$-subregularity of set-valued valued mappings between general Banach and Asplund spaces concentrating mainly on subdifferential mappings. In addition to deriving verifiable conditions that imply and are implied by these notions, we compare them (when appropriate) with the corresponding notions of metric regularity and provide several numerical examples illustrating the new phenomena.

Section~4 studies behavior of strong metric $q$-subregularity as $q\ge 1$ under parameter perturbations. The obtained results, being of their own interest, allow us to establish in Section~5 the convergence rate for Newton's and quasi-Newton methods of solving generalized equations depending on the order of the strong metric subregularity for the underlying set-valued mapping in the generalized equation under consideration. Section~6 presents concluding remarks and some directions of our future research.\vspace*{0.05in}

Throughout the paper we use standard notation of variational analysis and generalized differentiation. Recall that, given a set-valued mapping $F\colon X\tto X^*$ from the Banach space $X$ into its topological dual $X^*$ endowed with the weak$^*$ topology $w^*$, the symbol
\begin{eqnarray}\label{pk}
\disp\Limsup_{x\to\ox}F(x):=\big\{x^*\in X^*\big|\;\exists\,\mbox{seqs. }\,x_k\to\ox,\;x^*_k\st{w^*}{\to} x^*\;\mbox{with }x^*_k\in F(x_k),\;k\in\N\Big\}
\end{eqnarray}
signifies the {\em sequential Painlev\'e-Kuratowski outer limit} of $F$ as $x\to\ox$, where $\N:=\{1,2,\ldots\}$. Given a set $\O\subset\R^n$ and an extended-real-valued function $\ph\colon\R^n\to\oR:=(-\infty,\infty]$ finite at $\ox$, the symbols $x\st{\O}{\to}\ox$ and $x\st{\ph}{\to}\ox$
stand for $x\to\ox$ with $x\in\O$ and for $x\to\ox$ with $\ph(x)\to\ph(\ox)$, respectively. As usual, $\B(x,r)=\B_r(x)$ denotes the closed ball of the space in question centered at $x$ with radius $r>0$, while the symbols $\B$ and $\B^*$ signify the corresponding closed unit ball in the primal and dual spaces, respectively. Finally, given a mapping $g\colon X\to Y$ between Banach spaces that is locally Lipschitzian around $\ox$, we denote
$$
\lip g(\ox):=\disp\limsup_{x,u\to\ox}\frac{\|g(x)-g(u)\|}{\|x-u\|}.
$$

\section{Generalized Differentiation}\sce

In this section we present for the reader's convenience some basic tools of generalized differentiation widely employed in what follows. We refer to the books \cite{bz05,bm06,rw,s} for more details in both finite and infinite dimensions. Since the subdifferential and normal cone constructions are used below only in Asplund spaces, we confine ourselves to their definitions on this setting. Recall that a Banach space is {\em Asplund} if each of its separable subspace has a separable dual. This class of spaces is rather large including, in particular, every reflexive Banach space.

Given $\ph\colon X\to\oR$ with $\ox\in\dom\ph$, the {\em regular subdifferential} (known also as the presubdifferential and as the Fr\'echet or viscosity subdifferential) of $\ph$ at $\ox$ is defined by
\begin{eqnarray}\label{2.1}
\Hat\partial\ph(\ox):=\Big\{v\in\R^n\Big|\;\liminf_{x\to\ox}\frac{\ph(x)-\ph(\ox)-\la
v,x-\ox\ra}{\|x-\ox\|}\ge 0\Big\}.
\end{eqnarray}
It reduces to $\{\nabla\ph(\ox)\}$ if $\ph$ is Fr\'echet differentiable at $\ox$ and to the subdifferential of convex analysis if $\ph$ is convex, while the set  $\Hat\partial\ph(\ox)$ may often be empty for nonconvex and nonsmooth functions as, e.g., for $\ph(x)=-|x|$ at $\ox=0\in\R$. A serious disadvantage of (\ref{2.1}) is the failure of standard calculus rules required in variational analysis and its applications to optimization.

We come to the different picture while performing a limiting procedure/robust regularization over the mapping $x\mapsto\Hat\partial\ph(x)$ as $x\st{\ph}{\to}\ox$ vis the sequential outer limit (\ref{pk}), which gives us the (basic first-order) {\em subdifferential} of $\ph$ at $\ox$ defined by
\begin{eqnarray}\label{2.2}
\partial\ph(\ox):=\disp\Limsup_{x\st{\ph}{\to}\ox}\Hat\partial\ph(x)
\end{eqnarray}
and known also as the general, or limiting, or Mordukhovich subdifferential. In contrast to (\ref{2.1}), the set (\ref{2.2}) is often nonconvex (e.g., $\partial\ph(0)=\{-1,1\}$ for $\ph(x)=-|x|$) enjoying nevertheless comprehensive calculus based on variational/extremal principles of variational analysis.

\section{Necessary and Sufficient Conditions for $q$-Subregularity}\sce

Let us start this section with the basic definition of positive-order metric subregularity for arbitrary set-valued mappings between Banach spaces.

\begin{Definition}[\bf metric $q$-subregularity and strong $q$-subregularity]\label{subreg} Let $F:X\tto Y$ with $(\bar x,\bar y)\in\gph F$, and let $q>0$. We say that:

{\bf (i)} $F$ is {\sc metrically $q$-subregular} at $(\bar x,\bar y)$ if there are constants $\eta,\gg>0$ such that
\begin{equation}\label{3.1}
d\big(x;F^{-1}(\bar y)\big)\le\eta d^q\big(\bar y;F(x)\big)\;\mbox{ for all }\;x\in\B(\bar x;\gamma)
\end{equation}
The infimum over all constants/moduli $\eta>0$ for which {\rm(\ref{3.1})} holds with some $\gg>0$ is called the {\sc exact q-subregularity bound} of F at $(\bar x,\bar y)$ and is denoted by $\subreg^q F(\bar x,\bar y)$.

{\bf (ii)} $F$ is {\sc strongly metrically $q$-subregular} at $(\bar x,\bar y)$ if there are $\eta,\gamma>0$ such that
\begin{equation}\label{3.2}
\|x-\bar x\|\le\eta d^q\big(\bar y;F(x)\big)\;\mbox{ for all }\;x\in\B(\bar x;\gamma).
\end{equation}
The infimum over all $\eta>0$ for which {\rm(\ref{3.2})} holds with some $\gg>0$ is called the {\sc exact strong q-subregularity bound} of F at $(\bar x,\bar y)$ and is denoted by $\ssubreg^q F(\bar x,\bar y)$.
\end{Definition}

For brevity, in what follows we omit the adjective ``metric" for $q$-subregularity. It is easy to see from the definitions that the strong $q$-subregularity of $F$ at $(\ox,\oy)$ implies the corresponding $q$-subregularity of $F$. Furthermore, the validity of $\bar q$-subregularity (resp.\ strong $\bar q$-subregularity) of $F$ ar $(\ox,\oy)$ for the fixed number $\bar q>0$ ensures this property for any $0<q\le\bar q$.

Clearly, the larger $q$ in the above subregularity properties the better the corresponding estimate (error bound) in (\ref{3.1}) and (\ref{3.2}) is. The following simple one-dimensional example shows that it makes sense to consider the $q$-subregularity property of order $q>1$, in contrast to its metric regularity counterpart of such (higher) orders, even in the case of real functions.

\begin{Example}[\bf $q$-subregularity of higher order]\label{exmsq2} {\rm Consider the continuous function $f(x):=|x|^{\frac{1}{2}}$, $x\in\R$, which is not Lipschitz continuous around $\ox=0$. We have
$$
|x|\le|x^{\frac{1}{2}}|^q\;\mbox{for any $q\in(0,2]$ and all $x\in\B(0,1)$}.
$$
This shows that $f$ is strongly $q$-subregular at $(0,0)$ whenever $q\in(0,2]$.}
\end{Example}

The next example is more involved, being still one-dimensional, and reveals an interesting phenomenon: a set-valued mapping may not be metrically regular around the given point while it is metrically subregular at this point with some $q>1$. This example concerns in fact solution maps of the parametric generalized equations of type (\ref{ge}), which fails to have the metric regularity property in common situations; see Section~1.

\begin{Example}[\bf $q$-subregular but not metrically regular solution maps to parametric generalized equations]\label{exmsq} Consider the solution map
\begin{eqnarray}\label{ge1}
S(x)=\big\{y\in R^n\big|\;0\in f(x,y)+Q(y)\big\}
\end{eqnarray}
of the parametric generalized equation {\rm(\ref{ge})} with $f(x,y):=x$ and $Q\colon\R\tto\R$ given by
\begin{eqnarray*}
Q(y):=\left\{\begin{array}{rcl}
{[\frac{1}{2^{k+1}},\frac{1}{2^k}]}&\mbox{for}
&y=\frac{1}{(\sqrt[3]{2})^k},\\
\frac{1}{2^{k+1}}&\mbox{for}&y\in\Big(\frac{1}{(\sqrt[3]{2})^{k+1}},\frac{1}{(\sqrt[3]{2})^k}\Big),\\
0&\mbox{for}&y=0,\\
{[-\frac{1}{2^k},-\frac{1}{2^{k+1}}]}&\mbox{for}&y=-\frac{1}{(\sqrt[3]{2})^k},\\
-\frac{1}{2^{k+1}}&\mbox{for}&y\in\Big(-\frac{1}{(\sqrt[3]{2})^k},-\frac{1}{(\sqrt[3]{2})^{k+1}}\Big)
\end{array}\right.
\end{eqnarray*}
as depicted in Figure~{\rm 1}. Then $S$ is not metrically regular around $(0,0)$ while it is strongly $q$-subregular of any order
$q\in(0,2]$ at this point.
\end{Example}
\begin{figure}[ht]
\centering
\includegraphics[width=0.35\textwidth]{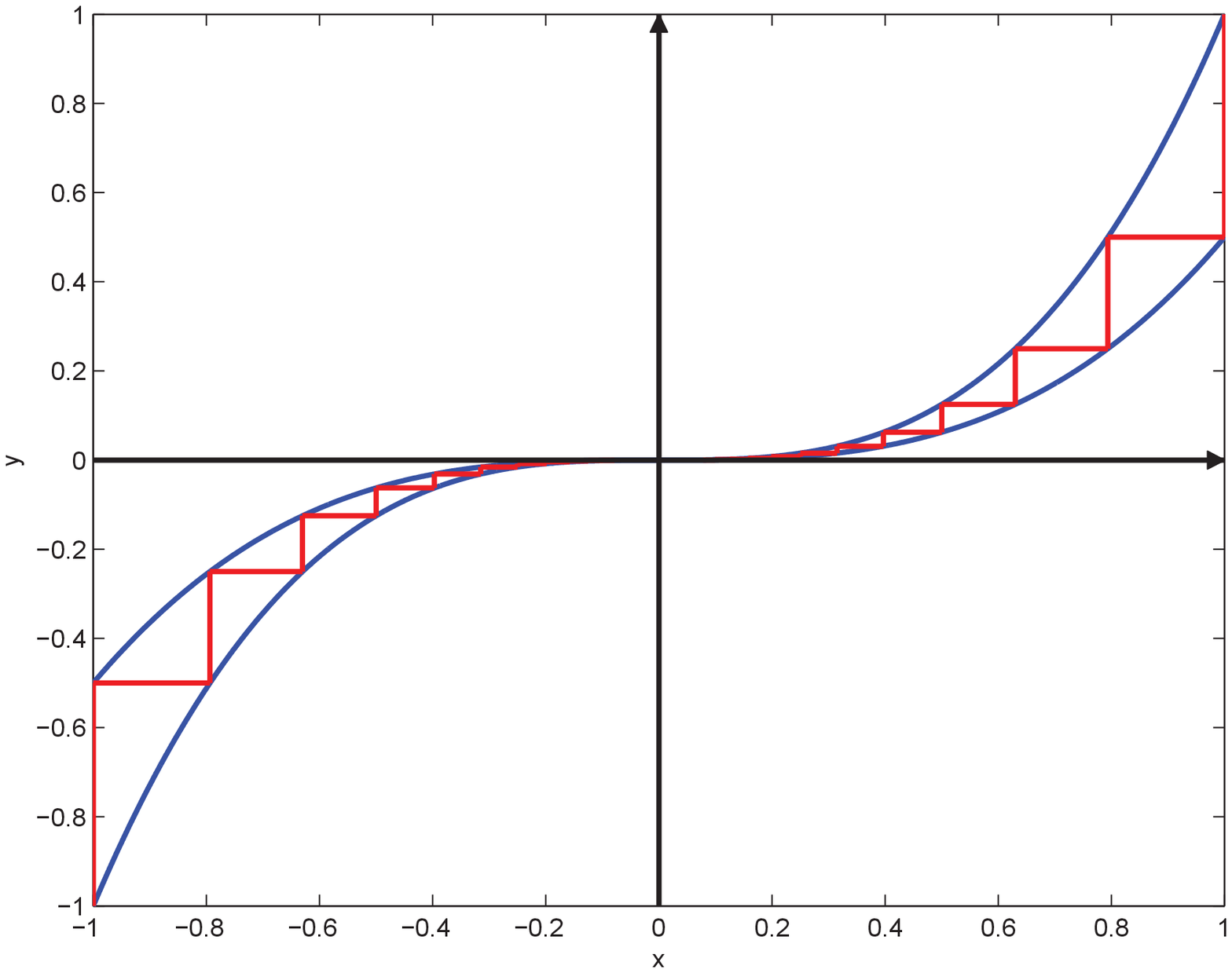}
\caption{Q(y)}
\label{fig:Q(y)}
\end{figure}

Indeed, the failure of metric regularity of $S$ in (\ref{ge1}) around $(0,0)$ follows from more general results of \cite{bm08}. Let us verify this directly for the mapping $S$ under consideration. Due to the form of $S$ in (\ref{ge1}) and the well-known equivalence between metric regularity of the given mapping and the Lipschitz-like/Aubin property of its inverse (see, e.g., \cite[Theorem~1.49]{bm06}), it suffices to show that $Q$ in (\ref{ge1}) is not Lipschitz-like around $(0,0)$. By \cite[Theorem~1.41]{bm06} this is equivalent to the fact that the scalar function
$$
\rho(x,y):=\text{dist}\big(x;Q(y)\big)=\inf\big\{\|x-v\|\big|\;v\in Q(y)\big\}
$$
is not locally Lipschitzian around $(0,0)$. To check the latter, we construct two sequences $\{(x_{1k},y_{1k})\}$ and $\{(x_{2k},y_{2k})\}$, which converge to $(0,0)$ when $k\rightarrow\infty$ as follows:
\begin{eqnarray*}
x_{1k}:=\frac{1}{2^{k-1}},\;y_{1k}:=\frac{1}{(\sqrt[3]{2})^{k-1}}-\alpha_k,\quad\text{where}\quad 0<\alpha_k<\min\Big\{\frac{1}{k2^k},\frac{1}{(\sqrt[3]{2})^{k-1}}-\frac{1}{(\sqrt[3]{2})^{k}}\Big\};
\end{eqnarray*}
$$
x_{2k}:=\frac{1}{2^{k-1}},\quad y_{2k}:=\frac{1}{(\sqrt[3]{2})^{k-1}},\quad k\in\N.
$$
Then we have the equalities
$$
\text{dist}(x_{2k};Q(y_{2k}))=\text{dist}\Big(x_{2k};{\Big[\frac{1}{2^k},\frac{1}{2^{k-1}}\Big]}\Big)=\text{dist}\Big(\frac{1}{2^{k-1}}, {\Big[\frac{1}{2^k},\frac{1}{2^{k-1}}\Big]}\Big)=0,
$$
$$
\text{dist}\big(x_{1k};Q(y_{1k})\big)=\text{dist}\Big(\frac{1}{2^{k-1}};\frac{1}{2^k}\Big)=\frac{1}{2^k},
$$
where the last one holds due to the estimates
$$
\frac{1}{(\sqrt[3]{2})^k}<y_{1k}<\frac{1}{(\sqrt[3]{2})^{k-1}},\quad k\in\N.
$$
Since $\|y_{1k}-y_{2k}\|=\alpha_k\le({k2^k})^{-1}$, we have
$$
\big\|d\big(x_{1k};Q(y_{1k})\big)-d\big(x_{2k};Q(y_{2k})\big)\big\|=\frac{1}{2^k}\ge k\alpha_k=k\|(x_{1k}-x_{2k},y_{1k}-y_{2k})\|,
$$
which indicates that $\rho(x,y)$ is not locally Lipschitzian around $(0,0)$ and yields therefore that the solution map $S$ from (\ref{ge1}) is not metrically regular around $(0,0)$.

Now we show that $S$ is $q$-subregular for $q=2$ and thus for any $q\in(0,2]$ at $(0,0)$. To proceed, take $\eta=\gamma=1$ and, given
$x\in B(0,\gamma)$, find $k_0$ so that $|x|\in{[\frac{1}{2^{k_0+1}},\frac{1}{2^{k_0}}]}$ and the corresponding value $|Q^{-1}(x)|$ belongs to ${\Big[\frac{1}{(\sqrt[3]{2})^{k+1}},\frac{1}{(\sqrt[3]{2})^{k-1}}\Big]}$. Notice that $S^{-1}(0)=\{0\}$, we have
\begin{eqnarray*}
d\big(x;S^{-1}(0)\big)=|x|&\le&\eta d^2\big(0;S(x)\big)=\eta\big(\inf\big\{\|y\|\big|\;y\in Q^{-1}(x)\big\}\big)^2,\quad x\in\B(0,\gamma),
\end{eqnarray*}
due to $\frac{1}{2^k}\le\frac{1}{(\sqrt[3]{2})^{2}(k+1)}$ for any $k\ge 2$, which verifies the $2$-subregularity of $S$ at $(0,0)$ and thus completes our justification in this example.\vspace*{0.05in}

It is worth mentioning that the solution map \eqref{ge1} in Example~\ref{exmsq} happens to be even {\em strongly} $2$-subregular at $(0,0)$. It follows from the arguments above since $S^{-1}(0)=\{0\}$ is a singleton.\vspace*{0.05in}

Next we derive {\em characterizations} of $q$-subregularity of any rate $q>0$ for the {\em subdifferential} mappings \eqref{2.2} generated by extended-real-valued lower semicontinuous (l.s.c.) functions on Banach (sufficient conditions) and Asplund (necessary conditions) spaces. For the case of subregularity ($q=1$) the obtained characterization reduces to \cite[Theorem~3.1]{dmn}. For convex functions on Banach spaces this case while concerning only local minimizers $x$ of $f$ has been independently characterized in \cite[Theorem~2.1]{ag} with a weaker modulus estimate; see more discussions in \cite{dmn} presented around Corollary~3.2 and also in \cite[Remark~2.2]{ag} for convex and nonconvex functions with $q=1$. In the general case of $q$-subregularity the formulation and proof of the theorem below are essentially more involved following the lines of the approach in \cite[Theorem~3.2]{mn} (for strong metric regularity) and of \cite[Theorem~3.1]{dmn} (for subregularity).

\begin{Theorem}[\bf characterization of $q$-subregularity of the basic subdifferential]\label{mqs} Let $f\colon X\rightarrow\overline {R}$ be l.s.c.\ around $\ox\in\dom f$ on a Banach space $X$, let $\bar x^*\in\partial f(\ox)$, and let $q$ be an arbitrary positive number. Consider the following two statements:

{\bf (i)} $\partial f$ is $q$-subregular at $(\ox,\ox^*)$ with modulus $\tilde\kappa$ and there exist numbers $\gamma>0$ and $r\in(0,q/\kappa)$ with $\kappa:=q\tilde{\kappa}^{\frac{1}{q}}$ such that
\begin{equation}\label{mqs1}
f(x)\ge f(\bar{x})+\langle{\bar x}^*,x-\bar{x}\rangle-\frac{qr}{1+q}d^{\frac{1+q}{q}}\big(x;(\partial f)^{-1}(\bar x^*)\big)\;\mbox{ for all }\;x\in \B(\bar x,\gamma).
\end{equation}

{\bf (ii)} There are two positive numbers $\alpha$ and $\eta$ such that
\begin{equation}\label{mqs2}
f(x)\ge f(\bar{x})+\langle{\bar x}^*,x-\bar{x}\rangle+\frac{q\alpha}{1+q}d^{\frac{1+q}{q}}\big(x;(\partial f)^{-1}(\bar x^*)\big)\;\mbox{ for all }\;x\in\B(\bar x,\eta).
\end{equation}
Then we have {\bf(ii)}$\Longrightarrow${\bf(i)} provided that there is $\beta\in(0,\al)$ with
\begin{equation}\label{mqs3}
f(u)\ge f(x)+\langle{x}^*,u-x\rangle-\frac{q\beta}{1+q}d^{\frac{1+q}{q}}\big(x,(\partial f)^{-1}(\bar x ^*)\big)
\end{equation}
whenever $(u,{\bar x}^*),(x,x^*)\in(\gph\partial f)\cap\B\big((\bar x,{\bar x}^*),\eta+(\frac{q\eta}{1+q})^{\frac{1}{q}}\big)$. Conversely, we have {\bf(i)}$\Longrightarrow${\bf(ii)} for any fixed $\al\in(0,q/\kappa)$ provided that the space $X$ is Asplund.
\end{Theorem}
{\bf Proof.} Let us first justify implication {\bf(ii)}$\Lra${\bf(i)} in the case of the general Banach space $X$ assuming condition \eqref{mqs3} with some $\beta\in(0,\alpha)$. Since {\rm(\ref{mqs2})} clearly yields {\rm(\ref{mqs1})}, we will arrive at {\bf (i)} by showing that there exists a number $\tilde\kappa>0$ such that
\begin{equation}\label{mqs7}
d\big(x;(\partial f)^{-1}(\bar x^*)\big)\le\tilde\kappa d^q\big(\bar x^*;\partial f(x)\big)\;\mbox{ for all }\;x\in\B\big(\bar x,\eta q/1+q\big).
\end{equation}
To proceed, fix $x\in\B(\bar x,\eta q/1+q)$. Pick any $u\in(\partial f)^{-1}(\bar x^*)$ with $\|x-u\|\le q^{-1}\|x-\bar x\|$ and get
$$
\|u-\bar x\|\le\|u-x\|+\|x-\bar x\|\le\eta.
$$
Then it follows from {\rm(\ref{mqs3})} the estimates
\begin{eqnarray}\label{mqs8}
f(\bar x)\ge f(u)+\langle{\bar x}^*,\bar x-u\rangle-\frac{\beta q}{1+q}d^{\frac{1+q}{q}}\big(u;(\partial f)^{-1}(\bar x^*)\big)\ge f(u)+\langle\bar x^*,\bar x-u\rangle,
\end{eqnarray}
which ensure in turn that for any such $u$ and $x^*\in\partial f(x)\cap\B(\bar x^*,[\eta q/(1+q)]^{\frac{1}{q}})$ we have
\begin{eqnarray*}
\langle x^*-\bar x^*,x-u\rangle&=&\langle x^*,x-u\rangle+\langle\bar x^*,u-\bar x\rangle-\langle\bar x^*,x-\bar x
\rangle\\&\ge&f(x)-f(u)-\frac{\beta q}{1+q}d^{\frac{q+1}{q}}\big(x;(\partial f)^{-1}(\bar x^*)\big)+f(u)-f(\bar x)-\langle\bar x^*,x-\bar x\rangle\\
&=&f(x)-f(\bar x)-\langle\bar x^*,x-\bar x\rangle-\frac{\beta q}{1+q}d^{\frac{q+1}{q}}\big(x;(\partial f)^{-1}(\bar x^*)\big)\\&\ge&\frac{q(\alpha-\beta)}{q+1}d^{\frac{q+1}{q}}\big(x;(\partial f)^{-1}(\bar x^*)\big),
\end{eqnarray*}
where the first inequality follows from {\rm(\ref{mqs3})} and {\rm(\ref{mqs8})} while the second one from {\rm(\ref{mqs2})}. Thus
$$
\|x^*-\bar x^*\|\cdot\|x-u\|\ge\frac{q(\alpha-\beta)}{q+1}d^{\frac{1+q}{q}}\big(x;(\partial f)^{-1}(\bar x^*)\big),
$$
which gives us the estimate
$$
\|x^*-\bar x^*\|d\big(x;(\partial f)^{-1}(\bar x^*)\big)\ge\frac{q(\alpha-\beta)}{q+1}d^{\frac{q+1}{q}}\big(x;(\partial f)^{-1}(\bar x^*)\big)
$$
due to the arbitrary choice of $u\in(\partial f)^{-1}(\bar x^*)$ with $\|x-u\|\le\frac{1}{q}\|x-\bar x\|\le\frac{\eta}{1+q}$. Hence
\begin{equation}\label{mqs9}
\|x^*-\bar x^*\|\ge\frac{q(\alpha-\beta)}{q+1}d^{\frac{1}{q}}\big(x;(\partial f)^{-1}(\bar x^*)\big)\;\mbox{ for all }\;x^*\in\partial f(x)\cap \B\big(\bar x^*,[\eta q/(1+q)]^{\frac{1}{q}}\big).
\end{equation}
If $d(\bar x^*;\partial f(x))\le(\eta q/1+q)^{\frac{1}{q}}$, we deduce from (\ref{mqs9}) that
$$
d\big(\bar x^*;\partial f(x)\big)\ge\frac{q(\alpha-\beta)}{q+1}d^{\frac{1}{q}}\big(x;(\partial f)^{-1}(\bar x^*)\big)\;\mbox{ for all }\;x\in\B(\bar x, \eta q/(1+q)\big),
$$
which justifies {\rm(\ref{mqs7})} with $\tilde\kappa:=[(1+q)/(q(\alpha-\beta)]^q$. In the remaining case of $d(\bar x^*;\partial f(x))>(\eta q/1+q)^{\frac{1}{q}}$ we obviously have the estimates
$$
d\big(x;(\partial f)^{-1}(\bar x^*)\big)\le\|x-\bar x\|\le\frac{\eta q}{1+q}<d^q\big(\bar x^*;\partial f(x)\big),
$$
which also justify {\rm(\ref{mqs7})} is and thus completes the proof of implication {\bf(ii)}$\Longrightarrow${\bf(i)}.

Next we verify the converse {\bf(i)}$\Lra${\bf(ii)} assuming that $X$ is Asplund. Arguing by contradiction, suppose that {\bf(i)} holds while property {\rm(\ref{mqs2})} is not satisfied whenever $\alpha,\eta>0$. Choose
$$
0<\frac{2^{\frac{1+q}{q}}}{1+q}<a<\infty
$$
and pick $\theta$ from the interval $\big(\frac{2^{-q}}{q(1+q)},\frac{1}{2}\big)$, which ensures that $\theta+\frac{1}{\theta^qa^q(1+q)^q}<1$. Now we claim that there exists a positive number $\nu$ satisfying
\begin{equation}\label{mqs4}
f(x)\ge f(\bar{x})+\langle\bar x^*,x-\bar{x}\rangle+\frac{q-(aq+1)\kappa r}{a(1+q)\kappa}d^{\frac{1+q}{q}}\big(x;(\partial f)^{-1}(\bar x^*)\big)\;\mbox{ for all }\;x\in\B(\bar x,\nu).
\end{equation}
Indeed, otherwise there is a sequence ${x_k}\rightarrow\bar{x}$ such that
\begin{equation*}
f(x_k)<f(\bar{x})+\langle{\bar x}^*,x_k-\bar{x}\rangle+\frac{q-(aq+1)\kappa r}{a(1+q)\kappa}d^{\frac{1+q}{q}}\big(x_k;(\partial f)^{-1}(\bar x^*)\big),\quad k\in\N.
\end{equation*}
This together with {\rm({\ref{mqs1}})} implies that all $x_k$ lie outside of $(\partial f)^{-1}(\bar x^*)$. Consequently we have
\begin{eqnarray}\label{mqs5}
\begin{array}{ll}
&\disp\inf_{x\in\B(\bar x,\gamma)}\Big\{f(x)-\langle{\bar x}^*,x-\bar{x}\rangle+\frac{q r}{1+q}d^{\frac{1+q}{q}}\big(x;(\partial f)^{-1}(\bar
x^*)\big)\Big\}\ge f(\bar{x})\\
&\disp>f(x_k)-\langle {\bar x}^*,x_k-\bar{x}\rangle+\frac{q r}{1+q}d^{\frac{1+q}{q}}\big(x_k;(\partial f)^{-1}(\bar x ^*)\big)-
\frac{q-\kappa r}{a(1+q)\kappa}d^{\frac{1+q}{q}}\big(x_k;(\partial f)^{-1}(\bar x ^*)\big).
\end{array}
\end{eqnarray}
Denote $\epsilon_k:=\big(\frac{q-\kappa r}{a(1+q)\kappa}\big)d^{\frac{1+q}{q}}\big(x_k;(\partial f)^{-1}(\bar x^*)\big)\dn 0$ as $k\to\infty$ and define the function
$$
g(x):=f(x)-\langle{\bar x}^*,x-\bar{x}\rangle+\frac{qr}{1+q}d^{\frac{1+q}{q}}\big(x;(\partial f)^{-1}(\bar x^*)\big),\quad x\in X.
$$
It follows from {\rm({\ref{mqs5}})} that $g(x_k)<\inf_{x\in\B(\bar x,\gamma)}g(x)+\epsilon_k$. Applying Ekeland's variational principle (see, e.g., \cite[Theorem~2.26]{bm06}) to the function $g+\delta_{\B(\bar x,\gamma)}$ with $\lambda_k:=\theta d\big(x_k;(\partial f)^{-1}(\bar x^*)\big)$ ensures the existence of a new sequence $\Hat{x}_k$ satisfying $\|\hat{x}_k-x_k\|\le\lambda_k$ and such that for each $k\in\N$ we have $\Hat{x}_k\in\text{int}\B(\bar x,\gamma)$ (due to $x_k\rightarrow\bar x$ and $\lambda_k\downarrow 0$ as $k\to\infty$) and that
$$
g(\hat{x}_k)<g(x)+\frac{\epsilon_k}{\lambda_k}\|x-\hat{x}_k\|\;\mbox{ for all }\;x\in\B(\bar x,\gamma).
$$
Employing the Fermat stationary rule in the above optimization problem and then using the subdifferential sum rule held in Asplund spaces by \cite[Theorem~3.41]{bm06}, we arrive at
\begin{eqnarray}0 &\in&\partial\Big(g(\cdot)+\frac{\epsilon_k}{\lambda_k}\|\cdot-\hat{x}\|\Big)(\hat{x}_k)\nonumber\\
&\subset&-{\bar x}^*+\partial f(\hat{x}_k)+\frac{\epsilon_k}{\lambda_k}\B^*+\frac{q r}{1+q}\partial d^{\frac{1+q}{q}}\big(\cdot;(\partial f)^{-1}(\bar x ^*)\big)(\hat{x}_k)\nonumber\\
&\subset&-{\bar x}^*+\partial f(\hat{x}_k)+\Big(r d^{\frac{1}{q}}\big(\hat{x}_k;(\partial f)^{-1}(\bar x^*)\big)+\frac{\epsilon_k}{\lambda_k}\Big) \B^*\nonumber.
\end{eqnarray}
Combining this with the metric q-subregularity property of $\partial f$ at $(\bar x,{\bar x}^*)$ ensures the estimates
$$
d^{\frac{1}{q}}\big(\hat{x}_k;(\partial f)^{-1}(\bar{x}^*)\big)\le\frac{\kappa}{q}d\big(\bar{x}^*;\partial f(\hat{x}_k)\big)
\le\frac{\kappa}{q}\Big(rd^{\frac{1}{q}}\big(\hat{x}_k;(\partial f)^{-1}(\bar x^*)\big)+\frac{\epsilon_k}{\lambda_k}\Big)
$$
for all $k\in\N$ sufficiently large. Hence for such numbers $k$  we get the inequality
$$
\Big(1-\frac{\kappa r}{q}\Big)d^{\frac{1}{q}}\big(\hat{x}_k;(\partial f)^{-1}(\bar x^*)\big)\le\frac{\kappa}{q}\frac{\epsilon_k}{\lambda_k}.
$$
This allows us to successively deduce that
\begin{eqnarray*}
\Big(1-\frac{\kappa r}{q}\Big)^q d\big(x_k;(\partial f)^{-1}(\bar x^*)\big)&\le&\Big(1-\frac{\kappa r}{q}\Big)^q\|\hat{x}_k-x_k\|+\Big(1-\frac{\kappa r}{q}\Big)^q d\big(\hat{x}_k;(\partial f)^{-1}(\bar x^*)\big)\\&\le& \Big(1-\frac{\kappa r}{q}\Big)^q\lambda_k+\Big(\frac{\kappa\epsilon_k}{q\lambda_k}\Big)^q\\&\le&\Big(1-\frac{\kappa r}{q}\Big)^q d\big(x_k;(\partial f)^{-1}(\bar x^*)\big)\Big(\theta+\frac{1}{a^q\theta^q(q+1)^q}\Big)\\&<&\Big(1-\frac{\kappa r}{q}\Big)^q d\big(x_k;(\partial f)^{-1}(\bar x^*)\big),
\end{eqnarray*}
where the last strict inequality follows from our choices of $a$ and $\theta$. Thus we arrive at the obvious contradiction, which justifiers our claim in {\rm({\ref{mqs4}})}. We conclude therefore that $q-(aq+1)\kappa r\le 0$.

Define now the real number
$$
r_1:=\frac{\frac{aq+1}{q}\kappa r-1}{a\kappa}\in\big[0,q/\kappa\big)
$$
and observe that inequality {\rm({\ref{mqs4}})} can be transformed into {\rm({\ref{mqs1}})} with replacing $r$ by $r_1$ and $\gamma$ by $\nu$, respectively. Consequently there is some real number $\nu_1$ such that
$$
f(x)\ge f(\bar{x})+\langle{\bar x}^*,x-\bar{x}\rangle+\frac{q-(aq+1)\kappa r_1}{a(1+q)\kappa}d^{\frac{1+q}{q}}\big(x;(\partial f)^{-1}(\bar x^*)\big)\;\mbox{ for all }\;x\in\B(\bar x,\nu_1).
$$
As before, we get the inequality $\kappa r_1>\frac{q}{aq+1}$, or equivalently $\kappa r>\frac{q}{aq+1}+\frac{aq^2}{(aq+1)^2}$. Defining
$$
r_2:=\frac{\frac{aq+1}{q}\kappa r_1-1}{a\kappa}\in\big[0,q/\kappa\big)
$$
and proceeding in the same way as above lead us to the inequality $\kappa r>\frac{q}{aq+1}+\frac{aq^2}{(aq+1)^2}+\frac{a^2q^3}{(aq+1)^3}$. Then we get by induction the progressively stronger bounds
$$
\kappa r>\frac{q}{aq+1}+\frac{aq^2}{(aq+1)^2}+\ldots+\frac{a^{k-1}q^k}{(aq+1)^k}=q\Big(1-\Big(\frac{aq}{aq+1}\Big)^k\Big)\;\mbox{ for all }\;k\in\N.
$$
Letting $k\rightarrow\infty$ gives us $\kappa r>q$, which is a contradiction. Therefore there exist real numbers $\alpha,\eta>0$ such that inequality {\rm({\ref{mqs2}})} is satisfied.

To justify {\bf(ii)}, we verify now that $\alpha$ may be chosen arbitrarily close to $q/\kappa$ while being smaller than this number. It suffices to consider the case of $\alpha<q/\kappa$. Take $\Tilde a,\tilde\theta>0$ such that
$$
\frac{q2^{\frac{1}{q}}}{\tilde{a}(1+q)}<\frac{1}{2},\quad\tilde{\theta}\in\Big(\frac{2^{\frac{1}{q}}}{\tilde{a}(1+q)},\frac{1}{2}\Big),\;\mbox{ and so }\;\tilde{\theta}+\Big(\frac{q}{\tilde{\theta}\tilde{a}(1+q)}\Big)^q<1.
$$
Given $\alpha,\eta>0$ for which \eqref{mqs2} holds, let us prove the existence $\mu>0$ such that
\begin{equation}\label{mqs6}
f(x)\ge f(\bar{x})+\langle{\bar x}^*,x-\bar{x}\rangle+\frac{q(q+(\tilde{a}-1)\alpha\kappa)}{\tilde{a}(1+q)\kappa}d^{\frac{1+q}{q}}\big(x;(\partial f)^{-1}(\bar x^*)\big)\;\mbox{ for all }\;x\in\B(\bar x,\mu).
\end{equation}
We just sketch the proof of (\ref{mqs6}) observing that it is similar to the proof of {\rm({\ref{mqs4}})} given above. Arguing by contradiction, find a sequence of $u_k\rightarrow\bar{x}$ so that
\begin{equation*}
f(u_k)<f(\bar{x})+\langle{\bar x}^*,u_k-\bar{x}\rangle+\frac{q(q+(\tilde{a}-1)\alpha\kappa)}{\tilde{a}(1+q)\kappa}d^{\frac{1+q}{q}}\big(u_k;(\partial f)^{-1}(\bar x ^*)\big),\quad k\in\N.
\end{equation*}
This gives us together with {\rm({\ref{mqs2}})} that
\begin{eqnarray*}
\begin{array}{ll}
&\disp\inf_{x\in\B(\bar x,\eta)}\Big\{f(x)-\langle{\bar x}^*,x-\bar{x}\rangle-\frac{q\alpha}{1+q}d^{\frac{1+q}{q}}\big(x;(\partial f)^{-1}(\bar x ^*)\big)\Big\}\ge f(\bar{x})\\
&\disp>f(u_k)-\langle{\bar x}^*,u_k-\bar{x}\rangle-\frac{q\alpha}{1+q}d^{\frac{1+q}{q}}\big(u_k;(\partial f)^{-1}(\bar x^*)\big)-\frac{q(q-\kappa\alpha)}{\tilde{a}(1+q)\kappa}d^{\frac{1+q}{q}}\big(u_k;(\partial f)^{-1}(\bar x^*)\big).
\end{array}
\end{eqnarray*}
Denote further $\nu_k:=\frac{q(q-\kappa\alpha)}{\tilde{a}(1+q)\kappa}d^{\frac{1+q}{q}}\big(u_k;(\partial f)^{-1}(\bar x^*)\big)\dn 0$ as $k\to\infty$ and consider the function
$$
h(x):=f(x)-\langle{\bar x}^*,x-\bar{x}\rangle-\frac{q\alpha}{1+q}d^{\frac{1+q}{q}}\big(x;(\partial f)^{-1}(\bar x^*)\big),\quad x\in X,
$$
for which we have $h(u_k)<\inf_{x\in\B(\bar x,\eta)}h(x)+\nu_k$ whenever $k\in\N$. Applying Ekeland's variational principle to the function $h+\delta_{\B(\bar x,\eta)}$ with $\rho_k:=\tilde{\theta}d(u_k;(\partial f)^{-1}(\bar x^*))$ ensures the existence of a new sequence $\{\hat{u}_k\}$ satisfying $\|\hat{u}_k-u_k\|\leq\rho_k$ and such that $\hat{u}_k\in\text{int}\,\B(\bar x,\eta)$ with
$$
h(\hat{u}_k)<h(x)+\frac{\nu_k}{\rho_k}\|x-\hat{u}_k\|\;\mbox{ for all }\;x\in\B(\bar x,\gamma),\;k\in\N.
$$
By the calculus rules as above we get the inclusions
\begin{equation*}
0\in\partial\Big(h(\cdot)+\frac{\nu_k}{\rho_k}\|\cdot-\hat{x}\|\Big)(\hat{u}_k)\subset-{\bar x}^*+\partial f(\hat{u}_k)+\Big(\alpha d^{\frac{1}{q}} \big(\hat{u}_k;(\partial f)^{-1}(\bar x^*)\big)+\frac{\nu_k}{\rho_k}\Big)\B^*,
\end{equation*}
which ensure together with the $q$-subregularity of $\partial f$ at $(\bar x,{\bar x}^*)$ that
$$
d^{\frac{1}{q}}\big(\hat{x}_k;(\partial f)^{-1}(\bar{x}^*)\big)\le\frac{\kappa}{q}d\big(\bar{x}^*,\partial f(\hat{x}_k)\big)
\le\frac{\kappa}{q}\Big(\alpha d^{\frac{1}{q}}\big(\hat{u}_k;(\partial f)^{-1}(\bar x^*)\big)+\frac{\nu_k}{\rho_k}\Big)
$$
for all $k\in\N$ sufficiently large. Thus for such $k$ we arrive at the estimate
$$
\Big(1-\frac{\kappa\alpha}{q})d^{\frac{1}{q}}\big(\hat{u}_k;(\partial f)^{-1}(\bar x^*)\big)\le\frac{\kappa}{q}\frac{\nu_k}{\rho_k}.
$$
This allows us to successively deduce that
\begin{eqnarray*}
\Big(1-\frac{\kappa\alpha}{q}\Big)^q d\big(u_k;(\partial f)^{-1}(\bar x^*)\big)&\le&
\Big(1-\frac{\kappa\alpha}{q}\Big)^q\|\hat{u}_k-u_k\|+\Big(1-\frac{\kappa\alpha}{q}\Big)^q d\big(\hat{u}_k;(\partial f)^{-1}(\bar x^*)\big)\\
&\le&\Big(1-\frac{\kappa\alpha}{q}\Big)^q\rho_k+\Big(\frac{\kappa\nu_k}{q\rho_k}\Big)^q\\
&\le&\Big(1-\frac{\kappa\alpha}{q}\Big)^q d\big(u_k;(\partial f)^{-1}(\bar x^*)\big)\Big(\tilde{\theta}+\Big(\frac{q}{\tilde{a}\tilde{\theta}(q+1)}\Big)^q\Big)\\
&<&\Big(1-\frac{\kappa\alpha}{q}\Big)^q d\big(u_k;(\partial f)^{-1}(\bar x^*)\big),
\end{eqnarray*}
which is a contradiction justifying the existence of $\mu>0$ such that {\rm({\ref{mqs6}})} holds.

In the last part of the proof we define the number
$$
\alpha_1:=\frac{(q+(\tilde{a}-1)\alpha\kappa)}{\tilde{a}\kappa}\in\big(0,q/\kappa\big)
$$
and observe that (\ref{mqs6}) yields {\rm({\ref{mqs2}})} with this number $\alpha_1$ and $\eta=\mu$ from (\ref{mqs6}). Then proceeding as above allows us to find $\mu_1>0$ such that
\begin{equation*}
f(x)\ge f(\bar{x})+\langle{\bar x}^*,x-\bar{x}\rangle+\frac{q(q+(\tilde{a}-1)\alpha_1\kappa)}{\tilde{a}(1+q)\kappa}d^{\frac{1+q}{q}}\big(x;(\partial f)^{-1}(\bar x^*)\big)\;\mbox{ for all }\;x\in\B(\bar x,\mu_1).
\end{equation*}
Define further $\alpha_2:=\frac{(q+(\tilde{a}-1)\alpha_1\kappa)}{\tilde{a}\kappa}\in(0,q/\kappa)$ and deduce from the above inequality that {\rm({\ref{mqs2}})} holds for $\alpha_2$ and $\eta=\mu_1$. By induction we find sequences $\{\alpha_k\}$ and $\{\mu_k\}$ satisfying
\begin{equation*}
f(x)\ge f(\bar{x})+\langle{\bar x}^*,x-\bar{x}\rangle+\frac{q}{1+q}\alpha_kd^{\frac{1+q}{q}}\big(x;(\partial f)^{-1}(\bar x^*)\big)\;\mbox{ for all }\;x\in B(\bar x,\mu_k)
\end{equation*}
with $\alpha_k:=\frac{(q+(\tilde{a}-1)\alpha_{k-1}\kappa)}{\tilde{a}\kappa}\in(0,q/\kappa)$ for $k\in\N$ and $\alpha_0=\alpha$. Letting finally $k\rightarrow\infty$ gives us $\alpha_k\rightarrow q/\kappa$ and thus completes proof of the theorem. $\h$\vspace*{0.05in}

Similarly to \cite[Corollaries~3.2, 3.3, 3.5]{dmn} given in the case of $q=1$, we can easily deduce from the obtained Theorem~\ref{mqs} its consequences for {\em local minimizers} as well as the characterization of {\em strong $q$-subregularity} in the case of arbitrary $q>0$. After completing and submitting this paper, we were informed by Xi Yin Zheng that implication (i)$\Longrightarrow$(ii) of Theorem~\ref{mqs} could be deduced from \cite[Theorem~4.1(i)]{zn} in the case when  $\ox$ is an isolated minimizer of $f$. If in addition to this the function $f$ is convex, then implication (ii)$\Longrightarrow$(i) could be deduced from \cite[Theorem~4.1(ii)]{zn} with the different proofs therein.

Now we present two examples illustrating the results of Theorem~\ref{mqs} and the assumptions made therein. The first example shows that the conditions of the theorem ensures the validity of 2-subregularity of the subdifferential mapping while metric regularity fails.

\begin{Example}[\bf 2-subregularity versus metric regularity of the subdifferential]\label{exsq} {\rm Consider the convex and continuous function $f\colon\R\to\R_+$ defined by
\begin{eqnarray*}
f(x):= \left\{\begin{array}{rcl}
-x&\mbox{for}
&x<-1,\\1&\mbox{for}&-1\le x\le 1,\\x&\mbox{for}&x>1.
\end{array}\right.
\end{eqnarray*}
It is not hard to calculate its subdifferential mapping as follows:
\begin{eqnarray*}
\partial f(x)=\left\{\begin{array}{rcl}
\{-1\}&\mbox{for}
&x<-1,\\{[-1,0]}&\mbox{for}&x=-1,\\\{0\}&\mbox{for}&-1<x<1,\\{[0,1]}&\mbox{for}&x=1,\\\{1\}&\mbox{for}&x>1.
\end{array}\right.
\end{eqnarray*}
Take $\bar x=0$, $\bar x^*=0\in\partial f(0)$ and show that the mapping $\partial f$ is 2-subregular (and hence $q$-subregular for any $0<q\le 2$) at $(0,0)$ while not metrically regular around this point. To verify 2-subregularity, it suffices to check by Theorem~\ref{mqs} that both conditions {\rm(\ref{mqs2})} and {\rm(\ref{mqs3})} of this theorem hold with $q=2$. To proceed, take $\alpha=1$ and $\eta=\frac{1}{2}$ in Theorem~\ref{mqs}(ii) and observe that $f(x)=f(0)=1$ for any $x\in\B(0,\eta)$. Then for any $q\in(0,2]$ we have
$$
f(x)\ge f(0)+\frac{1}{2}d^{1+\frac{1}{q}}\big(x;(\partial f)^{-1}(0)\big)
$$
due to $(\partial f)^{-1}(0)=(-1,1)$, and so condition {\rm(\ref{mqs2})} holds. To verify \eqref{mqs3}, take $\beta=0.5$ and get
$$
f(u)=f(x)=1\;\mbox{ for any }\;(u,\bar x^*),(x, x^*)\in\gph(\partial f)\cap\B\big((\bar x,\bar x^*),\eta+(\eta/2)^{\frac{1}{q}}\big),
$$
which justifies {\rm(\ref{mqs3})} and thus shows that $\partial f$ is 2-subregular at $(0,0)$ by Theorem~\ref{mqs}.

However, the validity of conditions \eqref{mqs2} and \eqref{mqs3} do not guarantee that the mapping $\partial f$ is metrically regular around $(0,0)$.
Indeed, we have for any $x,y\ne 0$ sufficiently close to $0$ that $\partial f(x)=\{0\}$ while the values of $(\partial f)^{-1}(y)$ are either $\{-1\}$ or $\{1\}$. Then it is easy to see by considering two sequences $\{x_k\}=\{k^{-1}\}$ and $\{y_k\}=\{(2k)^{-1}\}$ that there is no positive number $\kappa$, which ensures the distance estimate
$$
d\big(x_k;(\partial f)^{-1}(y_k)\big)\le\kappa d\big(y_k;\partial f(x_k)\big),\quad k\in\N,
$$
i.e., the subdifferential mapping $\partial f$ fails to be metrically regular around $(0,0)$.}
\end{Example}

The next example shows that condition {\rm(\ref{mqs3})} is {\em not necessary} for $q$-subregularity of $\partial f$ whenever $q>0$. In particular, this illustrates that (i) implies (ii) but does not imply \eqref{mqs3}.

\begin{Example}[\bf on assumptions and conclusions of Theorem~\ref{mqs}]\label{ex-mqs} {\rm
Define $f\colon\R\to\R_+$ by
\begin{eqnarray*}
f(x):=\left\{\begin{array}{rcl}
x^{\frac{1}{2}}&\mbox{for}&x>0,\\0&\mbox{for}&x=0,\\(-x)^{\frac{1}{2}}&\mbox{for}&x<0.
\end{array}\right.
\end{eqnarray*}
It follows immediately that $\partial f(x)=\{f'(x)\}$ at $x\ne 0$ with
\begin{eqnarray*}
f'(x)=\left\{\begin{array}{rcl}
\frac{1}{2\sqrt x}&\mbox{for}&x>0,\\-\frac{1}{2\sqrt{-x}}&\mbox{for}&x<0
\end{array}\right.
\end{eqnarray*}
and that ${[-1,1]}\subset\partial f(0)$. The latter implies that $\partial f^{-1}(0)=\{0\}$. Considering now $\bar x=0$ and $\bar x^*=0$, we claim that the mapping $\partial f$ is $q$-subregular at $(0,0)$ with any positive order $q$, which we fix in what follows. To proceed, take $\gg:=2^{-\frac{2q}{q+2}}$ and get
$$
|x|\le\big(1/2\big)^{\frac{2q}{2+q}}\Longleftrightarrow|x|\le\big(1/(2\sqrt{|x|})\big)^q\;\mbox{ for }\;x\in\B(0,\gg),
$$
which verifies the validity of the $q$-subregularity condition \eqref{3.1} for the mapping $\partial f$.

However, condition {\rm(\ref{mqs3})} fails here whenever $\beta,\eta>0$. To justify it, we argue by contradiction and suppose that {\rm(\ref{mqs3})} holds with some $\beta_0,\eta_0$. It is easy to see that the inclusion
$$
(u,0),(x,x^*)\in\gph(\partial f)\cap\B\Big((0,0),\eta_0+(\frac{q\eta_0}{1+q})^{\frac{1}{q}}\Big)
$$
yields $u=0$, which implies in turn that condition {\rm(\ref{mqs3})} reduces to
\begin{equation}\label{sc3}
0\ge f(x)+\la x^*,-x\ra-\frac{q\beta_0}{1+q}|x|^{\frac{q+1}{q}}
\end{equation}
Considering $x>0$ in {\rm(\ref{sc3})} gives us the inequality
$$
0\ge\frac{\sqrt x}{2}-\frac{q\beta_0}{1+q}|x|^{\frac{q+1}{q}},\;\mbox{ i.e., }\;\frac{2q\beta_0}{1+q}|x|^{\frac{q+1}{q}}\ge x^{\frac{1}{2}},
$$
which is a contradiction, since the latter inequality is obviously violated for $x$ sufficiently small.}
\end{Example}

\section{Strong Metric q-Subregularity under Perturbations}\sce

It has been well recognized in the literature that metric subregularity, in contrast to metric regularity, is not robust/stable with respect to perturbations of the initial data; see, e.g., \cite{dr09}. In this section we show that the situation is different for {\em strong} subregularity and higher-order $q$-subregularity ($q\ge 1$). Namely, it is proved below that such strong $q$-subregularity is {\em stable} with respect to appropriate perturbations of the initial set-valued mapping by Lipschitzian single-valued ones. Furthermore, we estimate the exact bound of strong $q$-subregularity moduli together with the radius of perturbations that keep strong $q$-subregularity of the perturbed mappings. These results are used in Section~5 for establishing the convergence rate for a class of quasi-Newton methods depending on the order $q$ of the assumed strong $q$-subregularity of the initial mapping. Unless otherwise stated, we have $q\ge 1$ in the rest of this section.

\begin{Theorem}[\bf strong $q$-subregularity under Lipschitzian perturbations]\label{per-lip} Let $F\colon X\tto Y$ be a set-valued mapping between Banach spaces with $(\ox,\oy)\in\gph F$, and let $g:X\to Y$ be a single-valued perturbation locally Lipschitzian around $\ox$. If there exist $\kappa,\lambda\in(0,\infty)$ such that
$$
\ssubreg^q F(\bar x,\bar y)<\kappa\;\mbox{ and }\;\lip g(\bar x)<\lambda<\big(\kappa^{\frac{1}{q}}\big)^{-1},
$$
then we have the modulus upper estimate
\begin{eqnarray}\label{e1}
\ssubreg^q(\tilde F+g)(\bar x,\bar y)\le\frac{\kappa}{(1-\lambda\kappa^{\frac{1}{q}})^q}\;\mbox{ with }\;\tilde F(x):=F(x)-g(\bar x).
\end{eqnarray}
Furthermore, if $\ssubreg^q F(\bar x,\bar y)>0$, we have the modulus relationship
\begin{eqnarray}\label{e2}
\ssubreg^q(\tilde{F}+g)(\bar x,\bar y)\le\frac{\ssubreg^q F(\bar x,\bar y)}{\big(1-(\ssubreg^qF(\bar x,\bar y)\big)^{\frac{1}{q}}\lip g(\bar x)^q}
\end{eqnarray}
whenever $(\ssubreg F(\bar x,\bar y))^{\frac{1}{q}}\lip g(\bar x)<1$.
\end{Theorem}
{\bf Proof.} Since $\ssubreg^q F(\bar x,\bar y)<\kappa$ and $\lip g(\bar x)<\lambda$, there exists $\gamma\in(0,1)$ such that
$$
\|x-\bar x\|\le\kappa d^q\big(\bar y;F(x)\big
)\;\text{ and }\;\|g(x)-g(\bar x)\|\le\lambda\|x-\bar x\|,\quad x\in\B(\bar x,\gamma).
$$
For such $x$, take $z\in\tilde{F}(x)+g(x)$ and find $y\in F(x)$ such that $z-y=g(x)-g(\bar x)$. It yields
\begin{eqnarray*}
\begin{array}{ll}
\|x-\bar x\|^{\frac{1}{q}}\le\kappa^{\frac{1}{q}}\|y-\bar y\|\le\kappa^{\frac{1}{q}} \big(\|\bar y-z\|+\|y-z\|\big)\\\\
\le\kappa^{\frac{1}{q}}\|\bar y-z\|+\kappa^{\frac{1}{q}} \lambda\|x-\bar x\|\le\kappa^{\frac{1}{q}} \|\bar y-z\|+\kappa^{\frac{1}{q}}\lambda\|x-\bar x\|^{\frac{1}{q}},
\end{array}
\end{eqnarray*}
where the last inequality holds due to $q\ge 1$ and $\gamma<1$. It gives us the estimate
$$
\|x-\bar x\|\le\frac{\kappa}{(1-\lambda\kappa^{\frac{1}{q}})^q}\|\bar y-z\|^q\;\mbox{ for all }\;x\in\B(\bar x,\gamma),\ z\in\tilde F(x)+g(x),$$
which implies both inequalities \eqref{e1}, \eqref{e2} and so completes the proof of the theorem.$\h$\vspace*{0.05in}

Next we derive two useful consequences of Theorem~\ref{per-lip} of their own interest. The first one concerns strictly differentiable (in particular, $C^1$-smooth) perturbations and is employed to establish convergence rates of the quasi-Newton methods considered in Section~5..

\begin{Corollary}[\bf strong $q$-subregularity under smooth perturbations]\label{per-sm} Let in the setting of Theorem~{\rm\ref{per-lip}} the perturbation $g$ is strictly differentiable at $\ox$, and let $\bar y\in F(\bar x)+g(\bar x)$. Then the mapping $x\mapsto F(x)+g(x)$ is strongly $q$-subregular at $(\bar x,\bar y)$ if and only if the mapping $G\colon x\mapsto g(\bar x)+\triangledown g(\bar x)(x-\bar x)+F(x)$ is strongly $q$-subregular at $(\bar x,\bar y)$ with the exact strong $q$-subregularity bound $\ssubreg^q(F+g)(\bar x,\bar y)$.
\end{Corollary}
{\bf Proof.} Pick $\kappa\in(0,\infty)$ such that $\ssubreg^q(F+g)(\bar x)<\kappa$. Define the mapping
$$
\tilde g(x):=\triangledown g(\bar x)(x-\bar x)+g(\bar x)-g(x)
$$
and observe that $\tilde g(\bar x)=0$ and $\lip\tilde{g}(\bar x)=0$ since $g$ is strictly differentiable at $\bar x$. Thus it follows from Theorem~\ref{per-lip} that the mapping $x \mapsto \widetilde{(F+g)}(x)+\tilde g(x)=F(x)+g(x)-\tilde g(\bar x)+\tilde g(x)=G(x)$ is also strongly $q$-subregular at $(\bar x,\bar y)$ with the exact bound not exceeding $\kappa$. The proof of the converse implication is similar with replacing $\tilde g$ by $-\tilde g$. $\h$\vspace*{0.05in}

The second consequence of Theorem~\ref{per-lip} provides a lower estimate of moduli of Lipschitzian perturbations, which fails strong $q$-subregularity of the original mapping.

\begin{Corollary}[\bf perturbation radius for failure of strong $q$-subregularity] In the setting of Theorem~{\rm\ref{per-lip}} we have the following estimate:
\begin{eqnarray}\label{pesq1}
\inf_{g\colon X\to Y}\big\{\lip g(\bar x)\big|\;\tilde F+g\;\text{\  is\  not \ strongly\ metrically \ q-subregular\ }\
at \ (\bar x,\bar y) \ \big\}\\\nonumber\ge\frac{1}{\big(\ssubreg^q F(\bar x,\bar y)\big)^{\frac{1}{q}}}.
\end{eqnarray}
\end{Corollary}
{\bf Proof.} We split the proof into considering the three possible cases in the theorem.

{\bf (i)} If $\ssubreg^q F(\bar x,\bar y)=\infty$, then the right-hand side of {\rm(\ref{pesq1})} is zero. Observing that for $g\equiv 0$ the mapping $\tilde F$ is not strong $q$-subregular, we conclude that the infimum in {\rm(\ref{pesq1})} the is also zero, and thus the inequality therein holds.

{\bf (ii)} If $\ssubreg^q F(\bar x,\bar y)=0$, then the right-hand side of {\rm(\ref{pesq1})} becomes $\infty$. For any $g: X \rightarrow Y$ with
$\lip g(\bar x)<\infty$ we deduce from Theorem~\ref{per-lip} that the mapping $\tilde F+g$ is strongly $q$-subregular as well. Hence the infimum in {\rm(\ref{pesq1})} is  also $\infty$, and thus the inequality holds therein.

{\bf(iii)} Consider the major case of $0<\ssubreg^qF(\bar x,\bar y)<\infty$ and suppose that {\rm(\ref{pesq1})} is violated. Then we find a mapping
$g:X\rightarrow Y$ locally Lipschitzian around $\ox$ such that
$$
\lip g(\bar x)\big(\ssubreg^qF(\bar x,\bar y)\big)^{\frac{1}{q}}<1
$$
and the mapping $\tilde F+g$ is not strongly $q$-subregular at $(\bar x,\bar y)$. This clearly contradicts Theorem~\ref{per-lip} and thus completes the proof. $\h$\vspace*{0.05in}

The concluding result of this section concerns strong $q$-subregularity of parameterized mappings being important, in particular, in the framework of Section~5.

\begin{Theorem}[\bf strong $q$-subregularity of parameterized mappings]\label{str} Let $F:X\tto Y$ be as above, and let $g\colon X\rightarrow Y$ be $C^1$-smooth around $\ox$. Assume that the mapping $G\colon x\mapsto g(\bar x)+\triangledown g(\bar x)(x-\bar x)+F(x)$ is strongly $q$-subregular at $(\bar x,\bar y)$ with $\bar y\in G(\bar x)$. Then for any $\lambda>\ssubreg^q G(\bar x,\bar y)$ there exists $\gamma>0$ such that the parameterized form of $G$ defined by
$$
x\mapsto G(u,x):=g(\bar x)+\triangledown g(u)(x-\bar x)+F(x)\;\mbox{ with }\;u\in\B(\bar x,\gamma)
$$
is strongly $q$-subregular at $(\bar x,\bar y)$ with modulus $\lm$, i.e., there is $\eta>0$ for which
$$
\|x-\bar x\|\le\lambda d^q\big(\bar y;G(u,x)\big)\;\mbox{ whenever }\;x\in\B(\bar x,\eta).
$$
\end{Theorem}
{\bf Proof.} Take $\lambda>\kappa>\ssubreg^q(G,\bar x,\bar y)$ and select $\mu>0$ so that
$$
\lambda>\frac{\kappa}{(1-\mu\kappa^{\frac{1}{q}})^q}\;\mbox{ and }\;\mu\kappa^{\frac{1}{q}}<1.
$$
By the assumed $C^1$ property of $g$, find $\gamma>0$ such that
$$
\|\triangledown g(u)-\triangledown g(\bar x)\|\le\mu\;\mbox{ for all }\;u\in\B(\bar x,\gamma).
$$
Fix further $u$ as above and define a new parameterized mapping $\tilde g:X\rightarrow Y$ by
$$
\tilde g(x):=g(u)+\triangledown g(u)(x-u)-g(\bar x)-\triangledown g(\bar x)(x-\bar x).
$$
Then we have $\tilde g(\bar x)=g(u)+\triangledown g(u)(\bar x-u)-g(\bar x)$, and hence
\begin{eqnarray*}
\|\tilde g(x)-\tilde g(x')\|&=&\|\triangledown g(u)(x-x')-\triangledown g(\bar x)(x-x')\|\\
&\le&\|\triangledown g(u)-\triangledown g(\bar x)\|\cdot\|x-x'\|\le\mu\|x-x'\|
\end{eqnarray*}
for any $x,x'\in\B(\bar x,\gamma)$. Thus $\lip\tilde g(\bar x)\le\mu$. Applying Theorem~\ref{per-lip} to the mappings $G$ and $\tilde g$ ensures that
the mapping
\begin{eqnarray*}
\begin{array}{ll}
x\mapsto\tilde G(x)+\tilde g(x)=G(x)-\tilde g(\bar x)+\tilde g(x)\\
g(\bar x)+\triangledown g(\bar x)(x-\bar x)+F(x)+(\triangledown g(u)-\triangledown g(\bar x))(x-\bar x)\\
=g(\bar x)+\triangledown g(u)(x-\bar x)+F(x)=G(u,x)
\end{array}
\end{eqnarray*}
is strongly $q$-subregular at $(\bar x,\bar y)$ with th exact bound not exceeded $\lambda$. This tells us that
for any $u\in\B(\bar x,\gamma)$ there is $\eta>0$ such that
$$
\|x-\bar x\|\le\lambda d^q\big(\bar y;G(u,x)\big)\;\mbox{ for all }\; x\in\B(\bar x,\eta),$$
which thus completes the proof of the theorem. $\h$

\section{Applications to Quasi-Newton Methods}\sce

In this section we discuss some applications of the results on strong $q$-subregularity under perturbations obtained in Section~4 to the convergence rate for a class of {\em quasi-Newton methods} to solve generalized equations given in the form
\begin{equation}\label{5.2}
0\in g(x)+F(x),
\end{equation}
where $g:X\rightarrow Y$ is a single-valued mapping while $F\colon X\tto Y$ is a set-valued mapping between Banach spaces. As in \cite{dmd}, we consider the following class of quasi-Newton methods to solve \eqref{5.2}:
\begin{equation}\label{5.1}
0\in g(x_k)+B_k(x_{k+1}-x_k)+F(x_{k+1}),\quad k=0,1,\ldots,
\end{equation}
where $B_k$ signify a sequence of linear and bounded operators acting from $X$ to $Y$. In the case of $B_k=\nabla g(x_k)$ algorithm \eqref{5.1} corresponds to {\em Newton's method} while particular choices of the operator sequence $\{B_k\}$ make it possible to include in this scheme various versions of quasi-Newton methods; see more discussions in \cite{dm,dmd}.

Let $\bar x$ be a solution to {\rm(\ref{5.2})}, and let $\{x_k\}$ be a sequence generated by {\rm(\ref{5.1})} that converges to $\bar x$. The classical Dennis-Mor\'e theorem \cite[Theorem~2.2]{dm} for $F\equiv 0$ establishes a certain characterization of the superlinear convergence  (called in \cite{dm} the $Q$-superlinear convergence, where $Q$ stands for ``quotient") of the quasi-Newton iterations
\begin{eqnarray}\label{super}
\lim_{k\to\infty}\frac{\|x_{k+1}-\ox\|}{\|x_k-\ox\|}=0
\end{eqnarray}
under the smoothness of $g\colon\R^n\to\R^n$ and nonsingularity of its Jacobian $\nabla g(\ox)$. Recently \cite[Theorem~3]{dmd}, Dontchev extended this result to the case of generalized equations \eqref{5.2} assuming that the mapping $x\mapsto g(\bar x)+\triangledown g(\bar x)(x-\bar x)+F(x)$ is strongly subregular at $(\bar x,0)$, which reduces to the nonsingularity of $\nabla g(\ox)$ in the setting of \cite{dm}.

The following theorem imposes the $q$-subregularity of $F$ at $(\ox,-g(\ox))$ as $q\ge 1$ and shows, by using the approach somewhat different from both papers \cite{dm,dmd} and based on the stability results of Section~4, that we have the {\em higher convergence rate}
\begin{eqnarray}\label{super-q}
\lim_{k\to\infty}\frac{\|x_{k+1}-\ox\|}{\|x_k-\ox\|^q}=0,\quad q\ge 1,
\end{eqnarray}
which reduces to the superlinear one in \eqref{super} for $q=1$.

\begin{Theorem}[\bf convergence rate for quasi-Newton iterations]\label{nmq} Let $\ox$ be a solution of the generalized equation \eqref{5.2}, where
$g\colon X\to Y$ be a mapping between Banach spaces that is $C^1$-smooth on some convex neighborhood $U$ of $\ox$. Given a starting point $x_0\in U$ and a sequence of linear and bounded operators $B_k\colon X\to Y$, consider the corresponding sequence $\{x_k\}$ generated by \eqref{5.1} such that $\{x_k\}\subset U$ and $x_k\to\ox$ as $k\to\infty$. Assume that the set-valued mapping $F\colon X\tto Y$ in \eqref{5.2} is strongly $q$-subregular at $(\bar x,-g(\bar x))$ with some $q\ge 1$ and that there exist positive numbers $\kappa,\lambda$ for which
\begin{eqnarray}\label{super1}
\ssubreg^q F\big(\bar x,-g(\bar x)\big)<\kappa\;\mbox{ and }\;\lip g(\bar x)<\lambda<\big(\kappa^{\frac{1}{q}}\big)^{-1}.
\end{eqnarray}
Suppose also that $x_k\ne\bar x$ for all $k\in\N$. Then we have the implication
\begin{eqnarray}\label{super2}
\Big[\disp\lim_{k\to\infty}\frac{\big\|\big(B_k-\triangledown g(\bar x)\big)(x_{k+1}-x_k)\big\|}{\|x_{k+1}-x_k\|}=0\Big]\Longrightarrow\eqref{super-q}.
\end{eqnarray}
\end{Theorem}
{\bf Proof.} By \eqref{super1} it follows from Theorem~\ref{per-lip} that
$$
\ssubreg^q(\tilde F+g)\big(\bar x,-g(\bar x)\big)\le\frac{\kappa}{(1-\lambda\kappa^{\frac{1}{q}})^q}\;\mbox{ with }\;\tilde F(x)=F(x)-g(\bar x).
$$
Then Corollary~\ref{per-sm} tells us that the mapping $\tilde G\colon x\mapsto g(\bar x)+\triangledown g(\bar x)(x-\bar x)+\tilde F(x)=\triangledown g(\bar x)(x-\bar x)+F(x)$ is strongly $q$-subregular at $(\bar x,-g(\bar x)$ with modulus $\mu\le\frac{\kappa}{(1-\lambda\kappa^{\frac{1}{q}})^q}$. To prove implication \eqref{super2} with $q\ge 1$, take any small $\epsilon>0$ and find a natural number $k_0$ sufficiently large so that for all $k\ge k_0$ we have the relationships
\begin{eqnarray*}
\|x_{k+1}-\bar x\|^{\frac{1}{q}}&\le&\mu^{\frac{1}{q}} d\big(-g(\bar x);\tilde G(x_{k+1})\big)=\mu^{\frac{1}{q}} d\big(-g(\bar x);\triangledown g(\bar x)(x_{k+1}-\bar x)+F(x_{k+1})\big)\\
&\le&\mu^{\frac{1}{q}}\|-g(\bar x)+g(x_k)+B_k(x_{k+1}-x_k)-\triangledown g(\bar x)(x_{k+1}-\bar x)\|\\
&=&\mu^{\frac{1}{q}}\|g(x_k)-g(\bar x)-\triangledown g(\bar x)(x_k-\bar x)+\big(B_k-\triangledown g(\bar x)\big)(x_{k+1}-x_k)\|\\
&\le&\mu^{\frac{1}{q}}\Big\|\int_0^1\triangledown g\big(\bar x+t(x_k-\bar x)\big)(x_k-\bar x)\,dt-\triangledown g(\bar x)(x_k-\bar x)\Big\|\\
&+&\mu^{\frac{1}{q}}\|B_k-\triangledown g(\bar x)\|\|x_{k+1}-x_k\|\\
&\le&\mu^{\frac{1}{q}}\Big\|\int_0^1 \big(\triangledown g\big(\bar x+t(x_k-\bar x)\big)-\triangledown g(\bar x)\big)(x_k-\bar x)dt\Big\|+\mu^{\frac{1}{q}}\epsilon\|x_{k+1}-x_k\|\\
&\le&\mu^{\frac{1}{q}}\epsilon\|x_k-\bar x\|+\mu^{\frac{1}{q}}\epsilon\|x_{k+1}-\bar x\|+\mu^{\frac{1}{q}}\epsilon\|x_k-\bar x\|
\\
&\le&2\mu^{\frac{1}{q}}\epsilon\|x_k-\bar x\|+\mu^{\frac{1}{q}}\epsilon\|x_{k+1}-\bar x\|^{\frac{1}{q}}.
\end{eqnarray*}
This gives us the upper estimate
$$
\|x_{k+1}-\bar x\|\le\frac{2^q\mu\epsilon^q}{(1-\mu^{\frac{1}{q}}\epsilon)^q}\|x_k-\bar x\|^q,
$$
which ensures the validity of \eqref{super-q} and thus completes the proof of the theorem. $\h$\vspace*{0.05in}

It worth mentioning that the inverse implication also holds in \eqref{super2}, which in fact follows from the proof of \cite[Theorem~3]{dmd} for $q=1$; cf.\ also \cite{dm}.\vspace*{0.05in}

As shown by the examples of Section~3, strong higher-order subregularity $(q>1)$ holds in natural situations when metric regularity fails. The following simple example, where $F(x)$ is a non-Lipschitzian function, illustrates the application of Theorem~\ref{nmq} in such settings.

\begin{Example}[\bf quasi-Newton method for non-Lipschitzian $2$-subregular equations]\label{nm2} {\rm Let $g(x):=x^2$ and $F(x):=|x|^{\frac{1}{2}}$, $x\in\R$. Then $\bar x=0$ is a solution to the generalized equation {\rm(\ref{5.2})}, which in this case reduces to a nonsmooth equation defined by a non-Lipschitzian function. As has been well recognized in the literature (see, e.g., \cite{fp,is14} and the references therein), the vast majority of the results on Newton-type methods for nonsmooth equations concerns Lipschitzian ones, while non-Lipschitzian settings are highly challenging. Based on Example~\ref{exmsq2}, we conclude that $F$ is strongly $2$-subregular at $(0,0)$. It is easy to check that the other conditions of Theorem~\ref{nmq} are also satisfied. Consider now the quasi-Newton method \eqref{5.1} with
$$
B_k:=\frac{\big(2^{\frac{(k+1)!}{2}}\big)^{-1}+\big(2^{2k!}\big)^{-1}}{\big(2^{k!}\big)^{-1}-\big(2^{(k+1)!}\big)^{-1}}.
$$
It is easy to verify that $\lim_{k\rightarrow\infty}|B_k-\triangledown g(\bar x)|=0$. Then for any starting point $x_0$ close to $\bar x$, it follows that algorithm \eqref{5.1} generates a sequence $\{x_k\}=\{(2^{k!})^{-1}\}$ converging to $\ox=0$ with the convergence rate that exceeds $q=2$.}
\end{Example}

\section{Concluding Remarks}

This paper demonstrates, for the first time in the literature, that the notions of metric $q$-subregularity and strong metric $q$-subregularity for sent-valued mappings can be useful not only in the cases when $q\in(0,1]$ but also for $q>1$. This significantly differs metric $q$-subregularity from metric $q$-regularity, which does not make sense when $q>1$. Besides various examples, characterizations of these notions and their sensitivity analysis, we provide applications to the convergence rate of quasi-Newton methods of solving generalized equations. It seems that it is only the beginning in the study and applications of these fruitful higher-order notions of variational analysis. We intend to develop more applications of higher-order metric subregularity and strong subregularity to various aspects of optimization; in particular, to convergence rates of other important algorithms in numerical optimization.

\end{document}